\documentclass[11pt,bezier]{article}
\usepackage{amsmath}
\usepackage{amsfonts,amsthm,amssymb}
\usepackage{amsfonts}
\usepackage{graphics}
\usepackage{amssymb}
\newtheorem{lem}{Lemma}[section]
\newtheorem{thm}[lem]{Theorem}

\theoremstyle{definition}

\begin{document}
\title{On the sizes of $(k,l)$-edge-maximal $r$-uniform hypergraphs
\footnote{The research is supported by NSFC (Nos. 11531011, 11771039, 11771443).}}
\author{Yingzhi Tian$^{a}$ \footnote{Corresponding author. E-mail: tianyzhxj@163.com (Y. Tian), hjlai@math.wvu.edu (H. Lai), mjx@xju.edu.cn (J. Meng), murong.xu009@gmail.com (M. Xu).}, Hong-Jian Lai$^{b}$, Jixiang Meng$^{a}$, Murong Xu$^{c}$, \\
{\small $^{a}$College of Mathematics and System Sciences, Xinjiang
University, Urumqi, Xinjiang 830046, PR China}\\
{\small $^{b}$Department of Mathematics, West Virginia University,
Morgantown, WV 26506, USA}\\
{\small $^{c}$Department of Mathematics, The Ohio State University, Columbus,
OH 43210, USA}}

\date{}

\maketitle

\noindent{\bf Abstract } Let $H=(V,E)$ be a hypergraph, where $V$ is a  set of vertices and $E$ is a set of non-empty subsets of $V$ called edges. If all edges of $H$ have the same cardinality $r$, then $H$ is a $r$-uniform hypergraph; if $E$ consists of all $r$-subsets of $V$, then $H$ is a complete $r$-uniform hypergraph, denoted by $K_n^r$, where $n=|V|$.  A $r$-uniform hypergraph $H=(V,E)$ is $(k,l)$-edge-maximal if every subhypergraph $H'$ of $H$ with $|V(H')|\geq l$ has edge-connectivity at most $k$, but for any edge $e\in E(K_n^r)\setminus E(H)$, $H+e$ contains at least one subhypergraph $H''$ with $|V(H'')|\geq l$ and edge-connectivity at least $k+1$. In this paper, we obtain the lower bounds and the upper bounds of the sizes of $(k,l)$-edge-maximal hypergraphs. Furthermore, we show that these bounds are best possible. Thus  prior results in [Y.Z. Tian, L.Q. Xu, H.-J. Lai, J.X. Meng, On the sizes of $k$-edge-maximal $r$-uniform hypergraphs, arXiv:1802.08843v3] are extended.


\section{Introduction}

In this paper, we consider finite simple graphs. For graph-theoretical terminologies and notation not defined here, we follow \cite{Bondy}.
For a graph $G$, we use $\kappa'(G)$ to denote the $edge$-$connectivity$ of $G$. The $complement$ of a graph $G$ is denoted by $G^c$. For $X\subseteq E(G^c)$, $G+X$ is the graph with vertex set $V(G)$ and edge set $E(G)\cup X$. We will use $G+e$ for $G+\{e\}$. The $floor$ of a real number $x$, denoted by $\lfloor x\rfloor$, is the greatest integer not larger than $x$; the $ceil$ of a real number $x$, denoted by $\lceil x\rceil$, is the least integer greater than or equal to $x$. For two integers $n$ and $k$, we define $(_k^n)=\frac{n!}{k!(n-k)!}$ when $k\leq n$ and $(_k^n)=0$ when $k>n$.

For generalizing a prior result of Mader \cite{Mader},  Boesch and McHugh \cite{Boesch} introduced the following definitions. For integers $k$ and $l$  with $l>k\geq2$, a graph $G$ with $n=|V(G)|\geq l$ is a $(k,l)$-$graph$ if $\kappa'(G')\leq k$ for any $G'\subseteq G$ with $|V(G')|\geq l$. A $(k,l)$-graph $G$ is $(k,l)$-$edge$-$maximal$ if, for any $e\in E(G^c)$, $G+e$ has a subgraph $G'$ with $|V(G')|\geq l$ and $\kappa'(G')\geq k+1$. $(k,k+1)$-edge-maximal graphs have been studied in [6,8,9], among others.

\begin{thm}
Let $k\geq2$ be an integer, and $G$ be a $(k,k+1)$-edge-maximal graph on $n>k+1$ vertices. Each of the following holds.

(i) (Mader \cite{Mader}) $|E(G)|\leq(n-k)k+(_2^k)$. Furthermore, this bound is best possible.

(ii) (Lai \cite{Lai}) $|E(G)|\geq (n-1)k-(_2^k)\lfloor\frac{n}{k+2}\rfloor$. Furthermore, this bound is best possible.
\end{thm}

In \cite{Boesch}, Boesch and McHugh extended Theorem 1($i$) to $(k,l)$-edge-maximal graphs.

\begin{thm}(Boesch and McHugh \cite{Boesch})
Let $G$ be a graph of order $n$ and let $n\geq l\geq k+1$. Let $p,q\geq0$ be integers such that $n=p(l-1)+q$ with $0\leq q<l-1$. If $G$ is a $(k,l)$-edge-maximal graph, then

\[
|E(G)|\leq
\left\{
\begin{array}{ll}
\frac{p(l-1)(l-2)}{2}+(p-1+q)k, & l-1>2k\ and\ q\leq2k,\\
\frac{p(l-1)(l-2)}{2}+pk+\frac{q(q-1)}{2},& l-1>2k\ and\ q>2k,\\
\frac{(l-1)(l-2)}{2}+(n-l+1)k, &l-1\leq 2k.
\end{array}
\right.
\]

Furthermore, these bounds are best possible.
\end{thm}

In \cite{LaiZhang}, Lai and Zhang extended Theorem 1($ii$) to $(k,l)$-edge-maximal graphs.

\begin{thm}(Lai and Zhang \cite{LaiZhang})
Let $G$ be a graph of order $n$ and let $n\geq l\geq k+3\geq5$. Let $p,q\geq0$ be integers such that $n=p(l-1)+q$ with $0\leq q<l-1$. If $G$ is a $(k,l)$-edge-maximal graph, then

\[
|E(G)|\geq
\left\{
\begin{array}{ll}
\frac{(l-1)(l-2)}{2}+(n-l+1)k, &l\leq n<2k+4,\\
(n-1)k-\lfloor\frac{n}{k+2}\rfloor\frac{(k+1)^2-3(k+1)}{2},& l\leq2k+4\leq n,\\
(n-2a+1)k+a(a-1)-\lfloor\frac{n-2a}{k+2}\rfloor\frac{(k+1)^2-3(k+1)}{2},& n\geq l=2a\geq2k+5,\\
(n-2b)k+b^2-\lfloor\frac{n-2b-1}{k+2}\rfloor\frac{(k+1)^2-3(k+1)}{2},& n\geq l=2b+1\geq2k+5.
\end{array}
\right.
\]

Furthermore, these bounds are best possible.
\end{thm}

Let $H=(V,E)$ be a $hypergraph$, where $V$ is a finite set and $E$ is a set of non-empty subsets of $V$, called $edges$. An edge of cardinality 2 is just a graph edge. For a vertex $u\in V$ and an edge $e\in E$, we say $u$ is $incident$ $with$ $e$ or $e$ is $incident$ $with$ $u$ if $u\in e$.
If all edges of $H$ have the same cardinality $r$, then $H$ is a $r$-$uniform$ $hypergraph$; if $E$ consists of all $r$-subsets of $V$, then $H$ is a $complete$ $r$-$uniform$ $hypergraph$, denoted by $K_n^r$, where $n=|V|$. For $n<r$, the complete $r$-uniform hypergraph $K_n^r$ is just the hypergraph with $n$ vertices and no edges. The $complement$ of a $r$-uniform hypergraph $H=(V,E)$, denoted by $H^c$, is the $r$-uniform hypergraph with vertex set $V$ and edge set consisting of all $r$-subsets of $V$ not in $E$. A hypergraph $H'=(V',E')$ is called a $subhypergraph$ of $H=(V,E)$, denoted by $H'\subseteq H$, if $V'\subseteq V$ and $E'\subseteq E$.
For $X\subseteq E(H^c)$, $H+X$ is the hypergraph with vertex set $V(H)$ and edge set $E(H)\cup X$; for $X'\subseteq E(H)$, $H-X'$ is the hypergraph with vertex set $V(H)$ and edge set $E(H)\setminus X'$. We use $H+e$ for $H+\{e\}$ and $H-e'$ for $H-\{e'\}$ when $e\in E(H^c)$ and $e'\in E(H)$.
For $Y\subseteq V(H)$, we use $H[Y]$ to denote the hypergraph $induced$ by $Y$, where $V(H[Y])=Y$ and $E(H[Y])=\{e\in E(H): e\subseteq Y\}$. $H-Y$ is the hypergraph induced by $V(H)\setminus Y$.

For a hypergraph $H=(V,E)$ and two disjoint vertex subsets $X, Y\subseteq V$, let $E_H[X,Y]$ be the set of edges intersecting both $X$ and $Y$ and $d_H(X,Y)=|E_H[X,Y]|$. We use $E_H(X)$ and $d_H(X)$ for $E_H[X,V\setminus X]$ and $d_H(X,V\setminus X)$, respectively. If $X=\{u\}$, we use $E_H(u)$ and $d_H(u)$ for $E_H(\{u\})$ and $d_H(\{u\})$, respectively. We call $d_H(u)$ the $degree$ of $u$ in $H$.  The $minimum$ $degree$ $\delta(H)$ of $H$ is defined as $min\{d_H(u): u\in V\}$; the $maximum$ $degree$ $\Delta(H)$ of $H$ is defined as $max\{d_H(u): u\in V\}$. When $\delta(H)=\Delta(H)=k$, we call $H$ $k$-$regular$.

For a nonempty proper vertex subset $X$ of a hypergraph $H$, we call $E_H(X)$ an $edge$-$cut$ of $H$. The $edge$-$connectivity$ $\kappa'(H)$ of a hypergraph $H$ is $min\{d_H(X):\O\neq X\subsetneqq V(H)\}$. By definition, $\kappa'(H)\leq \delta(H)$. We call a hypergraph $H$ $k$-$edge$-$connected$ if $\kappa'(H)\geq k$. A hypergraph is connected if it is 1-edge-connected. A maximal connected subhypergraph of $H$ is called a $component$ of $H$. A $r$-uniform hypergraph $H=(V,E)$ is $(k,l)$-$edge$-$maximal$ if every subhypergraph $H'$ of $H$ with $|V(H')|\geq l$ has edge-connectivity at most $k$, but for any edge $e\in E(H^c)$, $H+e$ contains at least one subhypergraph $H''$ with $|V(H'')|\geq l$ and edge-connectivity at least $k+1$. 
If $H$ is a $(k,l)$-edge-maximal $r$-uniform hypergraph with $n=|V(H)|<l$, then $H\cong K_n^r$.
For results on the connectivity of hypergraphs, see [1,4,5] for references.

In order to construct the complete $r$-uniform hypergraph with the maximum number  of vertices and degree at most $k$, we introduce the parameter $t=t(k,r)$, which is determined by $k$ and $r$.

\noindent{\bf Definition 1.} For two integers $k$ and $r$  with $k,r\geq2$, define $t=t(k,r)$ to be the largest integer such that $(^{t-1}_{r-1})\leq k$. That is, $t$ is the integer satisfying $(^{t-1}_{r-1})\leq k<(^{t}_{r-1})$.

In \cite{Tian}, the authors determined, for given integers $n$, $k$ and $r$, the extremal sizes of $(k,t)$-edge-maximal $r$-uniform hypergraphs on $n$ vertices.

\begin{thm} (Tian, Xu, Lai and Meng \cite{Tian})
Let $H$ be a $(k,t)$-edge-maximal $r$-uniform hypergraph such that $n\geq t$ and $k, r\geq2$, where $n=|V(H)|$ and $t=t(k,r)$.  Then each of the following holds.

(i) $|E(H)|\leq (^{t}_{r})+(n-t)k$. Furthermore, this bound is best possible.

(ii) $|E(H)|\geq (n-1)k -((t-1)k-(^{t}_{r}))\lfloor\frac{n}{t}\rfloor$. Furthermore, this bound is best possible.
\end{thm}

The main goal of this research is to extend these results in \cite{Tian}.
For given integers $n$, $k$ and $r$, the extremal sizes of a $(k,l)$-edge-maximal $r$-uniform hypergraph on $n$ vertices are determined, where  $l\geq t+1$. Section 2 below is devoted to the study of some properties of $(k,l)$-edge-maximal $r$-uniform hypergraphs. In section 3, we give the upper bounds of the sizes of $(k,l)$-edge-maximal $r$-uniform hypergraphs and illustrate that these bounds are best possible. We obtain the lower bounds of the sizes of $(k,l)$-edge-maximal $r$-uniform hypergraphs and show that these bounds are best possible in section 4.

\section{Properties of $(k,l)$-edge-maximal $r$-uniform hypergraphs}

Because of Theorem 1.4, we assume $l\geq t+1$ in this paper.

\begin{lem}
Let $H=(V,E)$ be a $(k,l)$-edge-maximal $r$-uniform hypergraph such that $n\geq l\geq t+1$ and $k,r\geq2$, where $n=|V(H)|$ and $t=t(k,r)$. Assume $X$ is a proper nonempty subset of $V(H)$ such that $\kappa'(H)=|E_H(X)|$. Then each of the following holds.

(i) $E_{H^c}(X)\neq\O$.

(ii) $\kappa'(H)=|E_H(X)|=k$.
\end{lem}

\noindent{\bf Proof.} Let $n_1=|X|$ and $n_2=|V(H)\setminus X|$. Then $n=n_1+n_2$. Since $H$ is $(k,l)$-edge-maximal, we have $\kappa'(H)\leq k$.

($i$) Assume $E_{H^c}(X)=\O$. Then  $E_{H}(X)$ consists of all $r$-subsets of $V(H)$ intersecting both $X$ and $V(H)\setminus X$. Thus
$$|E_{H}(X)|=\sum_{s=1}^{r-1}(_s^{n_1})(_{r-s}^{n_2})=
(_r^{n})-(_r^{n_1})-(_r^{n_2}). $$
Let $g(x)=(_r^{x})+(_r^{n-x})$. It is routine to verify that $g(x)$ is a decreasing function when $1\leq x\leq n/2$. If $min\{n_1,n_2\}\geq2$, then by $min\{n_1,n_2\}\leq n/2$, we have
$$\kappa'(H)=|E_{H}(X)|=
(_r^{n})-(_r^{n_1})-(_r^{n_2})\geq (_r^{n})-(_r^{2})-(_r^{n-2})>(_{r-1}^{n-1})\geq \delta(H), \eqno(1) $$
which contradicts to $\kappa'(H)\leq\delta(H)$. Now we assume $min\{n_1,n_2\}=1$. Then
$$\kappa'(H)=|E_{H}(X)|=
(_r^{n})-(_r^{n_1})-(_r^{n_2})= (_r^{n})-(_r^{1})-(_r^{n-1})=(_{r-1}^{n-1})\geq \delta(H), $$
which implies $\kappa'(H)=\delta(H)=(_{r-1}^{n-1})$ and so $H$ is a complete $r$-uniform hypergraph. Thus $\kappa'(H)=(^{n-1}_{r-1})\geq (^{l-1}_{r-1})\geq (^{t}_{r-1})>k$, contrary to $\kappa'(H)\leq k$. Therefore
$E_{H^c}(X)\neq\O$ holds.

($ii$)
By ($i$), we have $E_{H^c}(X)\neq\O$.
Pick an edge $e\in E_{H^c}(X)$. Since $H$ is $(k,l)$-edge-maximal, there is a subhypergraph $H'\subseteq H+e$ such that $|V(H')|\geq l$ and
$\kappa'(H')\geq k+1$. We have $e\in H'$ by $H$ is $(k,l)$-edge-maximal. It follows that $(E_H(X)\cup\{e\})\cap E(H')$ is an edge-cut of $H'$. Thus $|E_H(X)|+1\geq |E_H(X)\cup\{e\}|\geq \kappa'(H')\geq k+1$, implying $\kappa'(H)=|E_H(X)|\geq k$. By $\kappa'(H)\leq k$, we obtain
$\kappa'(H)=|E_H(X)|=k$.
$\Box$

\begin{lem}
Let $H=(V,E)$ be a $(k,l)$-edge-maximal $r$-uniform hypergraph such that $n\geq l\geq t+1$ and $k,r\geq2$, where $n=|V(H)|$ and $t=t(k,r)$. Assume $X$ is a proper nonempty subset of $V(H)$ such that $|E_H(X)|=k$. Then each of the following holds.

(i) If $|X|\leq r-1$, then $H[X]$ contains no edges in $E(H)$, and each edge of $E_H(X)$ contains $X$ as a subset.

(ii) If $r\leq|X|\leq l-1$, then $H[X]$ is a complete $r$-uniform hypergraph and $|X|\geq t$.

(iii) If $|X|\geq l$, then $H[X]$ is also a $(k,l)$-edge-maximal $r$-uniform hypergraph.
\end{lem}

\noindent{\bf Proof.}
($i$) Since $H$ is a $r$-uniform hypergraph, $H[X]$ contains no edges in $E[H]$ if $|X|\leq r-1$. By Lemma 2.1 and $|E_H(X)|=k$, we obtain $k=|E_H(X)|\geq\delta(H)\geq\kappa'(H)=k$, implying each edge of $E_H(X)$ contains $X$ as a subset.

($ii$) Assume $r\leq|X|\leq l-1$.  If $H[X]$ is not complete, then there is an edge $e\in E(H[X]^c)\subseteq E(H^c)$ and so $H+e$ has no subhypergraph $H'$ with $|V(H')|\geq l$ and $\kappa'(H')\geq k+1$, contrary to the assumption that $H$ is $(k,l)$-edge-maximal. Hence $H[X]$ must be complete.

On the contrary, assume $|X|<t$. Since $\delta(H)\geq\kappa'(H)=|E_H(X)|=k$ and $(^{t-1}_{r-1})\leq k<(^{t}_{r-1})$,  in order to ensure each vertex in $X$ has degree at least $k$ in $H$, we must have $|X|=t-1$ and $k=(^{t-1}_{r-1})$. Moreover, each vertex in $X$ is incident with exact $(^{t-2}_{r-2})$ edges in $E_H(X)$, and thus $d_H(u)=k$ for each $u\in X$. By Lemma 2.1 ($i$), there is an $e$ intersecting both $X$ and $V(H)\setminus X$ but $e\notin E_H(X)$. Since $|X|\geq r$, there is a vertex $w\in X$ such that $w$ is not incident with $e$. Then $d_{H+e}(w)=k$. This implies $w$ is not contained in a $(k+1)$-edge-connected subhypergraph of $H+e$. But then each vertex in $X\setminus\{w\}$ has at most degree $k$ in $(H+e)-w$, and thus each vertex  in $X\setminus\{w\}$ is not contained in a $(k+1)$-edge-connected subhypergraph of $H+e$. This implies that there is no $(k+1)$-edge-connected subhypergraph with at least $l$ vertices in $H+e$, a contradiction. Thus we have
$|X|\geq t$.

($iii$) Assume $|X|\geq l$. If $H[X]$ is complete, then $\kappa'(H[X])=\delta(H[X])=(^{|X|-1}_{r-1})\geq(^{l-1}_{r-1})
\geq(^{t}_{r-1})>k$, contrary to the definition of $(k,l)$-edge-maximal hypergraph. Thus $H[X]$ is not complete.
For any edge $e\in E(H[X]^c)\subseteq E(H^c)$, $H+e$ has a subhypergraph $H'$ with $|V(H')|\geq l$ and $\kappa'(H')\geq k+1$. Since  $|E_H(X)|=k$, we have $E_H(X)\cap E(H')=\O$. As $e\in E(H')\cap E(H[X]^c)$, we conclude that $H'$ is a subhypergraph of $H[X]+e$, and so $H[X]$ is a $(k,l)$-edge-maximal $r$-uniform hypergraph.
$\Box$

\section{The upper bounds of the sizes of $(k,l)$-edge-maximal $r$-uniform hypergraphs}

We first extend the definition of star-like-$(k,l)$ graphs in \cite{Boesch} to hypergraphs.

\noindent{\bf Definition 2.} Let $k,l,r$ be integers such that $k,r\geq2$ and $l\geq t+1$, where $t=t(k,r)$. Star-like-$(k,l)$ $r$-uniform hypergraphs are defined constructively as follows. Start with a complete $r$-uniform hypergraph $K_{l-1}^r$. Call it the $nucleus$. Attach a single vertex $K_1$ or a complete $r$-uniform hypergraph $K_i^r$ $(r\leq i\leq l-1)$ to this nucleus using $k$ edges joining $K_1$ or $K_i^r$ to the nucleus. Call this attached hypergraph a  $satellite$. Attach an arbitrary number of such satellites to the nucleus in the same manner. We call $r$-uniform hypergraphs constructed in this manner $star$-$like$-$(k,l)$ $r$-$uniform$ $hypergraphs$. We use $SH(k,l,r)$ to denote the collection of all star-like-$(k,l)$ $r$-uniform hypergraphs. See Figure 1 for example.

\begin{center}
\scalebox{0.7}{\includegraphics{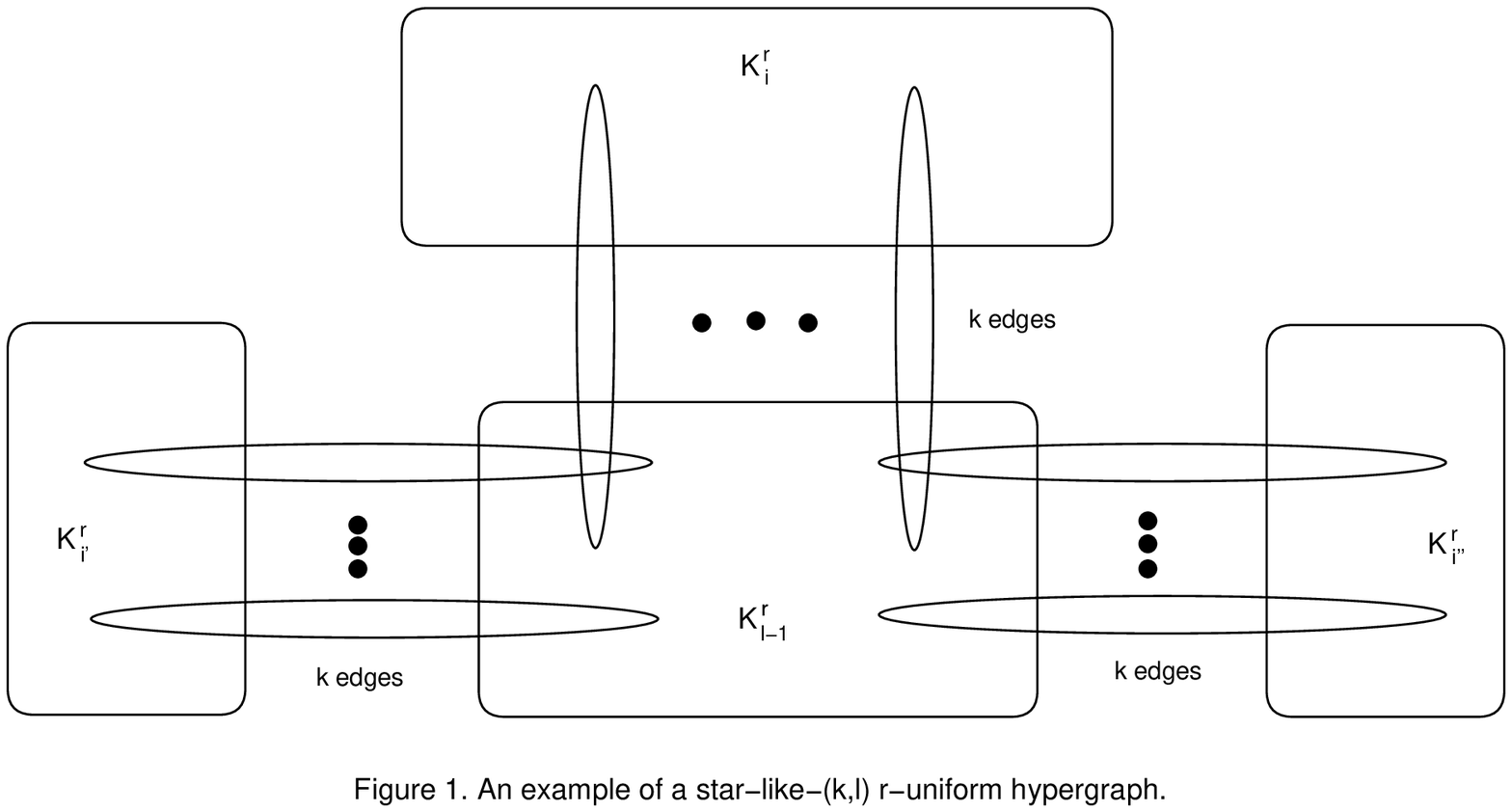}}
\end{center}

\noindent{\bf Definition 3.} For integers $k,r\geq2$, let $s=s(k,r)$ be the largest integer such that $k+(^s_r)\leq ks$.

\noindent{\bf Remark 1.} Since $k+(^t_r)=k+\frac{t}{r}(^{t-1}_{r-1})\leq k+\frac{t}{r}k\leq kt$, we have $t\leq s$, where $t=t(k,r)$ and $s=s(k,r)$.

\noindent{\bf Definition 4.} Let $n,k,l,r$ be integers such that $k,r\geq2$ and $n\geq l\geq t+1$, where $t=t(k,r)$. Let $p,q\geq0$ be integers such that $n=p(l-1)+q$ with $0\leq q<l-1$. We construct a class of star-like-$(k,l)$ $r$-uniform hypergraphs on $n$ vertices as follows.

($i$) If $l-1>s$ (where $s=s(k,r)$), then a star-like-$(k,l)$ $r$-uniform hypergraph consists of $p$ copies of $K_{l-1}^r$, one serving as the nucleus, the rest as satellites, together with addition satellites determined as follows: ($i$-$a$) if $q>s$, the single additional satellite is $K_{q}^r$; ($i$-$b$) if $q\leq s$, the additional satellites are $q$ copies of $K_1$, each attached to the nucleus by $k$ edges.

($ii$) If $l-1\leq s$ (where $s=s(k,r)$), then the nucleus is $K_{l-1}^r$. The satellites are $n-(l-1)$ copies of $K_1$, each attached to the nucleus by $k$ edges.

We denote the collection of all star-like-$(k,l)$ $r$-uniform hypergraphs on $n$ vertices constructed in Definition 4 by $MSH(n;k,l,r)$. Note that all hypergraphs in $MSH(n;k,l,r)$ have the same number of edges, denoted this number by $|E(MSH(n;k,l,r))|$ for brevity.  By definition, we have

\[
|E(MSH(n;k,l,r))|=
\left\{
\begin{array}{ll}
p(^{l-1}_r)+pk+(^{q}_r), & l-1>s\ and\ q>s,\\
p(^{l-1}_r)+(p-1+q)k,& l-1>s\ and\ q\leq s,\\
(^{l-1}_r)+(n-l+1)k, & l-1\leq s.
\end{array}
\right.
\]

The following theorem shows that $MSH(n;k,l,r)$ is a class of star-like-$(k,l)$ $r$-uniform hypergraphs with the maximum number of edges
among all star-like-$(k,l)$ $r$-uniform hypergraphs on $n$ vertices.

\begin{thm}
Let $n,k,l,r$ be integers such that $k,r\geq2$ and $n\geq l\geq t+1$, where $t=t(k,r)$. Let $p,q\geq0$ be integers such that $n=p(l-1)+q$ with $0\leq q<l-1$. For each star-like-$(k,l)$ $r$-uniform hypergraph $H$ on $n$ vertices,  we have $|E(H)|\leq |E(MSH(n;k,l,r))|$.
\end{thm}

\noindent{\bf Proof.} The idea of the proof is that we transform $H$ into one of the hypergraph in $MSH(n;k,l,r)$ by appropriate addition and deletion of edges. The edge transformations applied always yield an increase (not necessary net increase) in the number of edges. The proof uses two basic techniques: $splitting$ and $grouping$. The splitting operation replaces a single $K_i^r$ ($i\geq r$) by $i$ $K_1$-satellites. The grouping operations move vertices from smaller to larger satellites (with corresponding edge additions and deletions), or cluster a set of $K_1$-satellites into a single large satellite. We define the $satellite$ $spectrum$ of a star-like-$(k,l)$ hypergraph $H$ as $(S_1,S_r,\cdots,S_{l-1})$, where $S_i$, $i\in\{1,r,\cdots,l-1\}$, is the number of satellites of $H$ with $i$ vertices. We consider two cases.

\noindent{\bf Case 1.} $l-1>s$, where $s=s(k,r)$.

If $H$ contains a $K_i^r$-satellite, $r\leq i\leq s$, then perform the $splitting$ $operation$. That is, replace $K_i^r$-satellite (and the $k$ edges connecting it to the nucleus) by $i$ $K_1$-satellites (together with the $k$ edges that join each of them to the nucleus). The new hypergraph contains at least as many edges as $H$. The argument is as follows. The number of edges associated with the original $K_i^r$-satellite is $k+(^i_r)$. On the other hand, the satellites introduced by the splitting operation contribute $ik$ edges. Since $i\leq s$, we have  $k+(^i_r)\leq ik$ by Definition 3. Thus the hypergraph produced by the splitting operation contains at least  as many edges as $H$. If we repeat this process on any remaining $K_i^r$-satellites, $r\leq i\leq s$, we eventually obtain a transformed hypergraph with satellite spectrum satisfying $S_i=0$ for $r\leq i\leq s$.

If there are two satellites $K_i^r$ and $K_j^r$ satisfying $s<i\leq j<l-1$. We perform the following grouping operation (call $Grouping$ $operation$ $1$).  Since $i>s\geq t$ and $(^{i-1}_{r-1})\geq(^{t}_{r-1})>k$, we can assume that there is a vertex $u$ in $K_i^r$ not adjacent to the nucleus. Delete the edges connecting $u$ to $K_i^r$. Add edges from $u$ to $K_j^r$ such that $V(K_j^r)\cup\{u\}$ induces a complete $r$-uniform hypergraph. The operation is edge-increasing because while it removes $(^{i-1}_{r-1})$ edges, it adds $(^{j}_{r-1})$ edges. The $(K_i^r,K_j^r)$ pair of satellites become a $(K_{i-1}^r,K_{j+1}^r)$ pair. If $i-1\leq s$, then we apply the splitting operation to $K_{i-1}^r$. We repeat the grouping operation until at most one satellite remains in the range $s<i<l-1$. Now the satellite spectrum of the resulting hypergraph is simple. Except $S_1\geq0$, $S_{l-1}\geq0$, and some $S_{i_0}$ ($s<i_0<l-1$) may be 1 or 0, all other entries must be zero.

If $S_1$ and $S_{i_0}$ are positive, then we apply a second grouping operation (call $Grouping$ $operation$ 2). Select one $K_1$-satellite. Delete the $k$ edges connecting the $K_1$-satellite to the nucleus. Add $(^{i_0}_{r-1})$ edges connecting the $K_1$-satellite to $K_{i_0}^r$. This transformation is edge-increasing since $(^{i_0}_{r-1})>(^{s}_{r-1})\geq(^{t}_{r-1})>k$. We repeat this operation until the supply of $K_1$-satellites has been exhausted
or the original $K_{i_0}^r$ has been augmented to $K_{l-1}^r$-satellite.

If $S_1>s$, we apply one further grouping operation (call $Grouping$ $operation$ 3). Let $S_1=p_1(l-1) +q_1$, $0\leq q_1<l-1$. Replace the $S_1$ $K_1$-satellites by $p_1$ $K_{l-1}^r$-satellites, and by either $q_1$ $K_1$-satellites or 1 $K_{q_1}$-satellite, as $q_1\leq s$ or not. By definiton 3, this operation is edge increasing.

This completes the edge transformation of the original hypergraph in Case 1. For the resultant hypergraph, either $S_{i_0}=1$ for some $s<i_0<l-1$ and $S_1=0$, or all $S_i$ for $r\leq i\leq l-2$ are zero, which correspond to the case ($i$-$a$) or the case ($i$-$b$) in Definition 4.

\noindent{\bf Case 2.} $l-1\leq s$, where $s=s(k,r)$. We apply splitting operation to all satellites. By definition 3, the resultant hypergraph, which correspond to the case ($ii$) in Definition 4, is edge increasing.

This completes the proof of the theorem.
$\Box$

\begin{thm}
Let $n,k,l,r$ be integers such that $k,r\geq2$ and $n\geq l\geq t+1$, where $t=t(k,r)$. Let $p,q\geq0$ be integers such that $n=p(l-1)+q$ with $0\leq q<l-1$. If $H$ is a $(k,l)$-edge-maximal $r$-uniform hypergraph on $n$ vertices, then
\[
|E(H)|\leq |E(MSH(n;k,l,r))|=
\left\{
\begin{array}{ll}
p(^{l-1}_r)+pk+(^{q}_r), & l-1>s\ and\ q>s,\\
p(^{l-1}_r)+(p-1+q)k,& l-1>s\ and\ q\leq s,\\
(^{l-1}_r)+(n-l+1)k, & l-1\leq s,
\end{array}
\right.
\]
where $s=s(k,r)$.
\end{thm}

\noindent{\bf Proof.} By Theorem 3.1, we only need  to prove that there is a star-like-$(k,l)$ $r$-uniform hypergraph $H'$ on $n$ vertices such that $|E(H)|\leq|E(H')|$. The proof is by induction on $n$.

If $n=l$, let $H'$ be a star-like-$(k,l)$ $r$-uniform hypergraph with the nucleus $K_{l-1}^r$ and a single $K_1$-satellite. Since $\kappa'(K_{l}^r)=(^{l-1}_{r-1})$, we need to delete at least $(^{l-1}_{r-1})-k$ edges such that the remaining hypergraph have edge-connectivity at most $k$. Since $\kappa'(H)= k$ (by Lemma 2.1), we have $|E(K_{l}^r)|-|E(H)|\geq (^{l-1}_{r-1})-k$, implying $|E(H)|\leq|E(H')|$.

Now we assume $n>l$, and assume that for any $(k,l)$-edge-maximal $r$-uniform hypergraph with less than $n$ vertices, there is a star-like-$(k,l)$ $r$-uniform hypergraph having the same number of vertices and at least as many edges as the given hypergraph.

Let $F$ be a minimum edge-cut $H$. By Lemma 2.1, we have $|F|=k$. We consider two cases in the following.

\noindent{\bf Case 1.} There is a component, say $H_1$, of $H-F$ such that $|V(H_1)|=1$.

Let $H_2=H-V(H_1)$. Then $|V(H_2)|=n-1\geq l$. By Lemma 2.2 ($iii$), $H_2$ is $(k,l)$-edge-maximal. By induction assumption, there is a star-like-$(k,l)$ $r$-uniform hypergraph, say $H_2'$, such that $|V(H_2')|=|V(H_2)|$ and $|E(H_2')|\geq |E(H_2)|$. Let $H'$ be the star-like-$(k,l)$ $r$-uniform hypergraph obtained from $H_2'$ by adding a $K_1$-satellite. Since $|E(H)|=k+|E(H_2)|$, we have $|E(H)|\leq|E(H')|$.

\noindent{\bf Case 2.} Each component of $H-F$ has at least two vertices. Then, by lemma 2.2, each component of $H-F$ is either a complete $r$-uniform hypergraph with at least $t$ vertices, or a $(k,l)$-edge-maximal $r$-uniform hypergraph with at least $l$ vertices.

Let $H_1$ be a component of $H-F$ and $H_2=H-V(H_1)$. Assume $n_1=|V(H_1)|$ and $n_2=|V(H_2)|$. Then $n_1+n_2=n$.

\noindent{\bf Subcase 2.1.} $n_1\geq l$ and $n_2\geq l$.

By Lemma 2.2, both $H_1$ and $H_2$ are $(k,l)$-edge-maximal $r$-uniform hypergraphs.  By induction assumption,  there are two star-like-$(k,l)$ $r$-uniform hypergraphs, say $H_1'$ and $H_2'$, such that $|V(H_i')|=|V(H_i)|$ and $|E(H_i')|\geq |E(H_i)|$ for $i=1,2$. Let $H'$ be a star-like-$(k,l)$ $r$-uniform hypergraph obtained from $H_1'$ and $H_2'$ by moving the satellites of $H_1'$ to $H_2'$ and changing the nucleus of $H_1'$ to be a satellite of $H_2'$. Then $|E(H')|=|E(H_1')|+|E(H_2')|+k$, and thus $|E(H)|\leq|E(H')|$ holds.

\noindent{\bf Subcase 2.2.} $n_1\geq l$ and $n_2< l$.

By Lemma 2.2, $H_1$ is a $(k,l)$-edge-maximal $r$-uniform hypergraphs and $H_2$ is a complete $r$-uniform hypergraph with at least $t$ vertices.  By induction assumption,  there is  star-like-$(k,l)$ $r$-uniform hypergraphs, say $H_1'$, such that $|V(H_1')|=|V(H_1)|$ and $|E(H_1')|\geq |E(H_1)|$. Let $H'$ be a star-like-$(k,l)$ $r$-uniform hypergraph obtained from $H_1'$ by adding  $H_2$ to be a satellite of $H_1$. Then $|E(H')|=|E(H_1')|+|E(H_2)|+k$, and thus $|E(H)|\leq|E(H')|$ holds.

\noindent{\bf Subcase 2.3.} $n_1< l$ and $n_2< l$.

By Lemma 2.2, both $H_1$ and $H_2$ are complete $r$-uniform hypergraphs with at least $t$ vertices.  Since  $|E(H)|=(^{n_1}_r)+(^{n_2}_r)+k\leq (^{l-1}_r)+(^{n-l+1}_r)+k$, we obtain that each star-like-$(k,l)$ $r$-uniform hypergraph on $n$ vertices having at least as many edges as $H$ in this case. Therefore, the proof of this case follows.

This completes the proof of Theorem 3.2.
$\Box$

In the following theorem, we will show that each hypergraph in  $MSH(n;k,l,r)$ is $(k,l)$-edge-maximal. So the upper bounds given in Theorem 3.2 are best possible.

\begin{thm}
Let $n,k,l,r$ be integers such that $k,r\geq2$ and $n\geq l\geq t+1$, where $t=t(k,r)$.  If $H\in MSH(n;k,l,r)$, then $H$ is $(k,l)$-edge-maximal.
\end{thm}

\noindent{\bf Proof.} If $l-1=t$, then all satellites of $H$ are $K_1$-satellites by $s\geq t$, where $s=s(k,r)$.  By Lemma 3.1 in \cite{Tian}, $H$ is $(k,l)$-edge-maximal. Thus, in the following, we assume $l-1>t$.

By definition, there is no subhypergraph $H'$ of $H$ such that $|V(H')|\geq l$ and $\kappa'(H')>k$.  We will prove the theorem by induction on $n$. If $n=l$,
then $H$ is a star-like-$(k,l)$ $r$-uniform hypergraph with the nucleus $K^r_{l-1}$ and a $K_1$-satellite. Since $\kappa'(K^r_{l-1})=(^{l-2}_{r-1})\geq(^{t}_{r-1})>k$, for any $e\in E(H^c)$, we have $\kappa'(H+e)> k$. Thus $H$ is $(k,l)$-edge-maximal.

Now suppose $n>l$. We assume that each hypergraph in $MSH(n';k,l,r)$, where $n'<n$, is $(k,l)$-edge-maximal. In the following, we will show that each $H$ in $MSH(n;k,l,r)$ is also $(k,l)$-edge-maximal.

By contradiction, assume that there is an edge $e\in E(H^c)$ such that $H+e$ contains no subhypergraph $H'$ satisfying $|V(H')|\geq l$ and $\kappa'(H')>k$.
Let $F$ be an edge-cut in $H+e$ with cardinality at most $k$. Since $\kappa'(K^r_{l-1})=(^{l-2}_{r-1})\geq(^{t}_{r-1})>k$ and $\kappa'(K^r_{q})=(^{q-1}_{r-1})\geq(^{s}_{r-1})\geq(^{t}_{r-1})>k$ when $q>s$, we obtain that $F$ is exact the edge-cut joining some satellite and the nucleus. Thus there is a component, say $H_1$, of $H-F$, such that $H_1$ is the hypergraph obtained from $H$ by deleting one satellite and $e\in H_1^c$.
By induction assumption, $H_1+e$ contains a subhypergraph $H_1'$ such that $|V(H_1')|\geq l$ and $\kappa'(H_1')>k$. But $H_1'$ is also a subhypergraph of $H+e$, a contradiction.
$\Box$

\section{The Lower bounds of the sizes of $(k,l)$-edge-maximal $r$-uniform hypergraphs}

The following lemma will be needed in proving the main result in this section.

\begin{lem}
Let $n,a,k,r$ be integers such that $k,r\geq2$ and $n\geq a\geq t$, where $t=t(k,r)$.  We have the following two inequalities.

($i$) $(^n_r)\geq(n-1)k -((t-1)k-(^{t}_{r}))\lfloor\frac{n}{t}\rfloor$.

($ii$) $(^n_r)\geq(n-a)k+(^a_r)-((t-1)k-(^{t}_{r}))\lfloor\frac{n-a}{t}\rfloor$.
\end{lem}

\noindent{\bf Proof.} Since $|E(K^r_n)|=(^n_r)$, the lemma will hold if we can construct two $r$-uniform hypergraphs $H$ and $H'$ on $n$ vertices such that $|E(H)|=(n-1)k -((t-1)k-(^{t}_{r}))\lfloor\frac{n}{t}\rfloor$ and $|E(H')|=(n-a)k+(^a_r)-((t-1)k-(^{t}_{r}))\lfloor\frac{n-a}{t}\rfloor$.

Let $H$ be a $r$-uniform star-like hypergraph with the nucleus $K^r_t$, $\lfloor\frac{n}{t}\rfloor-1$ $K^r_t$-satellites and $n-t\lfloor\frac{n}{t}\rfloor$ $K_1$-satellites, adding $k$ edges joining each satellite to the nucleus.
It is routine to count that $|E(H)|=(n-1)k -((t-1)k-(^{t}_{r}))\lfloor\frac{n}{t}\rfloor$.

Let $H'$ be a $r$-uniform star-like hypergraph with the nucleus $K^r_a$, $\lfloor\frac{n-a}{t}\rfloor$ $K^r_t$-satellites and $n-a-t\lfloor\frac{n-a}{t}\rfloor$ $K_1$-satellites, adding $k$ edges joining each satellite to the nucleus.
Then  $|E(H')|=(n-a)k+(^a_r)-((t-1)k-(^{t}_{r}))\lfloor\frac{n-a}{t}\rfloor$.
$\Box$

It is routine to verify the following lemma.

\begin{lem}
For given integers $n$ and $r$ with $n\geq r\geq2$, the function $g(x)=(^x_r)+(^{n-x}_r)$ is decreasing in the range $1\leq x\leq n/2$.
\end{lem}

\begin{thm}
Let $n,k,l,r$ be integers such that $k,r\geq2$ and $n\geq l\geq t+1$, where $t=t(k,r)$. If $H$ is a $(k,l)$-edge-maximal $r$-uniform hypergraph on $n$ vertices, then
\[
|E(H)|\geq
\left\{
\begin{array}{ll}
(^{l-1}_r)+(n-l+1)k, &l\leq n<2t,\\
(n-1)k -((t-1)k-(^{t}_{r}))\lfloor\frac{n}{t}\rfloor, & l\leq2t\leq n,\\
(n-2a+1)k+2(^{a}_r)-((t-1)k-(^{t}_{r}))\lfloor\frac{n-2a}{t}\rfloor,& n\geq l=2a\geq2t+1,\\
(n-2b)k+(^{b}_r)+(^{b+1}_r)-((t-1)k-(^{t}_{r}))\lfloor\frac{n-2b-1}{t}\rfloor,& n\geq l=2b+1\geq2t+1.
\end{array}
\right.
\]

\end{thm}

\noindent{\bf Proof.} Let $F$ be a minimum edge-cut of $H$. By Lemma 2.1, $|F|=k$. Assume $H_1$ is a minimum component of $H-F$ and $H_2=H-V(H_1)$.
Let $n_1=|V(H_1)|$ and $n_2=|V(H_2)|$. Then $n=n_1+n_2$ and $n_1\leq n_2$.

($i$) For $l\leq n<2t$, we have $n_1<t$, and then $n_1=1$ by Lemma 2.2.
If $n=l$, then, by Lemma 2.2,  $H_2$ is a complete $r$-uniform hypergraph on $l-1$ vertices.
Thus $|E(H)|=k+(^{l-1}_r)=(^{l-1}_r)+(n-l+1)k$. If $n>l$, by Lemma 2.2, $H_2$ is $(k,l)$-edge-maximal. By induction on $n$, assume  $|E(H_2)|\geq(^{l-1}_r)+(n_2-l+1)k$. Therefore, $|E(H)|=|F|+|E(H_2)|\geq k+(^{l-1}_r)+(n_2-l+1)k=(^{l-1}_r)+(n-l+1)k$.

($ii$) We now assume $l\leq2t\leq n$. We shall prove this case by induction on $n$.

If $n=2t$, then either $n_1=n_2=t$, or $n_1=1$ and $n_2=n-1$ by Lemma 2.2. When $n_1=n_2=t$, then $H_1$ and $H_2$ are complete, and thus $|E(H)|=(^{t}_r)+(^{t}_r)+k=(n-1)k -((t-1)k-(^{t}_{r}))\lfloor\frac{n}{t}\rfloor$. Assume $n_1=1$ and $n_2=n-1$.
If $n=l$, then, by Lemma 2.2,  $H_2$ is a complete $r$-uniform hypergraph.
By Lemma 4.2, $(^{n-1}_r)=(^{1}_r)+(^{n-1}_r)\geq (^{t}_r)+(^{t}_r)$. 
Thus $|E(H)|=k+(^{n-1}_r)\geq k+2(^{t}_r)=(n-1)k -((t-1)k-(^{t}_{r}))\lfloor\frac{n}{t}\rfloor$. If $n>l$, by Lemma 2.2, $H_2$ is $(k,l)$-edge-maximal. Assume, by induction hypothesis, $|E(H_2)|\geq(n_2-1)k -((t-1)k-(^{t}_{r}))\lfloor\frac{n_2}{t}\rfloor$. Therefore,  $|E(H)|=|E(H_2)|+k\geq (n-1)k -((t-1)k-(^{t}_{r}))\lfloor\frac{n-1}{t}\rfloor\geq (n-1)k -((t-1)k-(^{t}_{r}))\lfloor\frac{n}{t}\rfloor,$ the last inequality holds because $(t-1)k-(^{t}_{r})\geq (t-1)(^{t-1}_{r-1})-\frac{t}{r}(^{t-1}_{r-1})\geq0$.

Assume that $n>2t$. If $n_1=1$, then $n_2=n-1\geq2t$. By induction assumption,
$|E(H_2)|\geq(n_2-1)k -((t-1)k-(^{t}_{r}))\lfloor\frac{n_2}{t}\rfloor$.
Thus $|E(H)|=|E(H_2)|+k\geq (n-1)k -((t-1)k-(^{t}_{r}))\lfloor\frac{n-1}{t}\rfloor\geq (n-1)k -((t-1)k-(^{t}_{r}))\lfloor\frac{n}{t}\rfloor$. So we assume $n_1\geq t$.

\noindent{\bf Claim.} $|E(H_i)|\geq(n_i-1)k -((t-1)k-(^{t}_{r}))\lfloor\frac{n_i}{t}\rfloor$ for $i\in \{1,2\}$.

If $|V(H_i)|\geq l$, then by induction assumption, we have  $|E(H_i)|\geq(n_i-1)k -((t-1)k-(^{t}_{r}))\lfloor\frac{n_i}{t}\rfloor$. If $t\leq|V(H_i)|\leq l-1$, then by Lemma 2.2, $H_i$ is a complete $r$-uniform hypergrpah. Thus, by Lemma 4.1 ($i$), $|E(H_i)|=(^{n_i}_r)\geq(n_i-1)k -((t-1)k-(^{t}_{r}))\lfloor\frac{n_i}{t}\rfloor$.

By this claim, we have

\ \ \ \ $|E(H)|=|E(H_1)|+|E(H_2)|+k$

\ \ \ \ \ \ \ \ \ \ \ \ \ \ $\geq(n_1-1)k -((t-1)k-(^{t}_{r}))\lfloor\frac{n_1}{t}\rfloor+(n_2-1)k -((t-1)k-(^{t}_{r}))\lfloor\frac{n_2}{t}\rfloor+k$

\ \ \ \ \ \ \ \ \ \ \ \ \ \ $=(n-1)k-((t-1)k-(^{t}_{r}))(\lfloor\frac{n_1}{t}\rfloor+
\lfloor\frac{n_2}{t}\rfloor)$

\ \ \ \ \ \ \ \ \ \ \ \ \ \ $\geq(n-1)k-((t-1)k-(^{t}_{r}))\lfloor\frac{n_1+n_2}{t}\rfloor$ (By $(t-1)k-(^{t}_{r})\geq0$)

\ \ \ \ \ \ \ \ \ \ \ \ \ \ $=(n-1)k-((t-1)k-(^{t}_{r}))\lfloor\frac{n}{t}\rfloor$.

($iii$) We then assume $n\geq l=2a\geq2t+1$. By induction on $n$, we will prove this case.

If $n=l$, then either $n_1=1, n_2=n-1$ or $t\leq n_1\leq n_2\leq n-t$. When $n_1=1$ and $n_2=n-1=l-1$, then $H_2$ is complete by Lemma 2.2. Since $(^{l-1}_r)=(^{1}_r)+(^{l-1}_r)=(^{1}_r)+(^{2a-1}_r)\geq (^{a}_r)+(^{a}_r)$ (by Lemma 4.2), we have  $|E(H)|=|E(H_2)|+k=k+(^{n-1}_r)\geq k+2(^{a}_r)=(n-2a+1)k+2(^{a}_r)-((t-1)k-
(^{t}_{r}))\lfloor\frac{n-2a}{t}\rfloor$.
When $t\leq n_1\leq n_2\leq l-t$, by Lemma 2.2, both $H_1$ and $H_2$ are complete.  Since $n_1\leq n/2=a$, we have $(^{n_1}_r)+(^{n_2}_r)\geq (^{a}_r)+(^{a}_r)$ by Lemma 4.2. 
Thus $|E(H)|=|E(H_1)|+|E(H_2)|+k=(^{n_1}_r)+(^{n_2}_r)+k\geq k+2(^{a}_r)=(n-2a+1)k+2(^{a}_r)-((t-1)k-
(^{t}_{r}))\lfloor\frac{n-2a}{t}\rfloor$.

Thus we assume $n>l$. If $n_1=1$, then $n_2=n-1\geq l$. By induction assumption, $|E(H_2)|\geq (n_2-2a+1)k+2(^{a}_r)-((t-1)k-
(^{t}_{r}))\lfloor\frac{n_2-2a}{t}\rfloor$. Thus
 $|E(H)|=|E(H_2)|+k\geq (n_2-2a+1)k+2(^{a}_r)-((t-1)k-
(^{t}_{r}))\lfloor\frac{n_2-2a}{t}\rfloor+k\geq (n-2a+1)k+2(^{a}_r)-((t-1)k-
(^{t}_{r}))\lfloor\frac{n-2a}{t}\rfloor$. So we assume $n_1\geq t$ and consider three cases in the following.

\noindent{\bf Case 1.} $n_1\geq l$.

By induction assumption, we have $|E(H_i)|\geq (n_i-2a+1)k+2(^{a}_r)-((t-1)k-
(^{t}_{r}))\lfloor\frac{n_i-2a}{t}\rfloor$ for $i=1,2$. By setting $n$ to be $a$ in Lemma 4.1 ($i$), we have $(^{a}_r)\geq(a-1)k-((t-1)k-(^{t}_{r}))\lfloor\frac{a}{t}\rfloor$.
Thus

\ \ \ \ $|E(H)|=|E(H_1)|+|E(H_2)|+k$

\ \ \ \ \ \ \ \ \ \ \ \ \ \ $\geq(n_1-2a+1)k+2(^{a}_r)-((t-1)k-
(^{t}_{r}))\lfloor\frac{n_1-2a}{t}\rfloor$

\ \ \ \ \ \ \ \ \ \ \ \ \ \
$+(n_2-2a+1)k+2(^{a}_r)-((t-1)k-
(^{t}_{r}))\lfloor\frac{n_2-2a}{t}\rfloor+k$

\ \ \ \ \ \ \ \ \ \ \ \ \ \ $\geq(n_1-2a+1)k+2(^{a}_r)-((t-1)k-
(^{t}_{r}))\lfloor\frac{n_1-2a}{t}\rfloor$

\ \ \ \ \ \ \ \ \ \ \ \ \ \
$+(n_2-2a+1)k+2((a-1)k-((t-1)k-(^{t}_{r}))\lfloor\frac{a}{t}\rfloor)-((t-1)k-
(^{t}_{r}))\lfloor\frac{n_2-2a}{t}\rfloor+k$

\ \ \ \ \ \ \ \ \ \ \ \ \ \ $\geq(n-2a+1)k+2(^{a}_r)-((t-1)k-
(^{t}_{r}))\lfloor\frac{n-2a}{t}\rfloor$.

\noindent{\bf Case 2.} $n_1<l$ and $n_2\geq l$.

By $n_1<l$ and $n_2\geq l$, we have  $H_1$ is complete and $H_1$ is $k$-edge-maximal by Lemma 2.2. By induction assumption,
$|E(H_2)|\geq (n_2-2a+1)k+2(^{a}_r)-((t-1)k-
(^{t}_{r}))\lfloor\frac{n_2-2a}{t}\rfloor$. Setting $n$ to be $n_1$ in Lemma 4.1 ($i$), we have $(^{n_1}_r)\geq(n_1-1)k-((t-1)k-(^{t}_{r}))\lfloor\frac{n_1}{t}\rfloor$ by $n_1\geq t$. Thus

\ \ \ \ $|E(H)|=|E(H_1)|+|E(H_2)|+k$

\ \ \ \ \ \ \ \ \ \ \ \ \ \ $\geq(^{n_1}_r)+(n_2-2a+1)k+2(^{a}_r)-((t-1)k-
(^{t}_{r}))\lfloor\frac{n_2-2a}{t}\rfloor+k$

\ \ \ \ \ \ \ \ \ \ \ \ \ \ $\geq(n_1-1)k-((t-1)k-(^{t}_{r}))\lfloor\frac{n_1}{t}\rfloor$

\ \ \ \ \ \ \ \ \ \ \ \ \ \
$+(n_2-2a+1)k+2(^{a}_r)-((t-1)k-
(^{t}_{r}))\lfloor\frac{n_2-2a}{t}\rfloor+k$

\ \ \ \ \ \ \ \ \ \ \ \ \ \ $\geq(n-2a+1)k+2(^{a}_r)-((t-1)k-
(^{t}_{r}))\lfloor\frac{n-2a}{t}\rfloor$.

\noindent{\bf Case 3.} $n_1\leq n_2<l$.

By $n_1\leq n_2<l$, we obtain that both $H_1$ and $H_2$ are complete by Lemma 2.2. If $n_1\geq a$, then by setting $n$ to be $n_i$ in Lemma 4.1 ($ii$), we have  $(^{n_i}_{r})\geq(n_i-a)k+(^{a}_r)-((t-1)k-
(^{t}_{r}))\lfloor\frac{n_i-a}{t}\rfloor$ for $i=1,2$. Thus

\ \ \ \ $|E(H)|=|E(H_1)|+|E(H_2)|+k$

\ \ \ \ \ \ \ \ \ \ \ \ \ \ $=(^{n_1}_{r})+(^{n_2}_{r})+k$

\ \ \ \ \ \ \ \ \ \ \ \ \ \ $\geq(n_1-a)k+(^{a}_r)-((t-1)k-
(^{t}_{r}))\lfloor\frac{n_1-a}{t}\rfloor$

\ \ \ \ \ \ \ \ \ \ \ \ \ \
$+(n_2-a)k+(^{a}_r)-((t-1)k-
(^{t}_{r}))\lfloor\frac{n_2-a}{t}\rfloor+k$

\ \ \ \ \ \ \ \ \ \ \ \ \ \
$\geq(n-2a+1)k+2(^{a}_r)-((t-1)k-
(^{t}_{r}))\lfloor\frac{n-2a}{t}\rfloor$.

If $n_1\leq a$, then $n_2=n-n_1\geq 2a-n_1\geq t$. By setting $n$ to be $n_2$ and $a$ to be $2a-n_1$ in Lemma 4.1 ($ii$), we have $(^{n_2}_{r})\geq(n_1+n_2-2a)k+(^{2a-n_1}_r)-((t-1)k-
(^{t}_{r}))\lfloor\frac{n_1+n_2-2a}{t}\rfloor$.  Together with  $(^{n_1}_r)+(^{2a-n_1}_r)\geq 2(^{a}_r)$ by Lemma 4.2, we have

\ \ \ \ $|E(H)|=|E(H_1)|+|E(H_2)|+k$

\ \ \ \ \ \ \ \ \ \ \ \ \ \ $=(^{n_1}_{r})+(^{n_2}_{r})+k$

\ \ \ \ \ \ \ \ \ \ \ \ \ \ $\geq(^{n_1}_{r})+(n_1+n_2-2a)k+(^{2a-n_1}_r)-((t-1)k-
(^{t}_{r}))\lfloor\frac{n_1+n_2-2a}{t}\rfloor+k$

\ \ \ \ \ \ \ \ \ \ \ \ \ \
$\geq(n-2a+1)k+2(^{a}_r)-((t-1)k-
(^{t}_{r}))\lfloor\frac{n-2a}{t}\rfloor$.

Therefore, the proof for ($iii$) of Theorem 4.3 is complete.

($iv$) The proof for ($iv$) of Theorem 4.3 is similar to that ($iii$) of Theorem 4.3, thus we omit the proof here.
$\Box$

\noindent{\bf Remark 2.} For $l\leq n<2t$, the star-like-($k,l$) hypergraphs in $MSH(n;k,l,r)$ show that the bound given in Theorem 4.3 ($i$) is best possible. Hypergraphs constructed in Definition 3 \cite{Tian} illustrate that the bound given in Theorem 4.3 ($ii$) is best possible. The following example will show that bounds given in Theorem 4.3 ($iii$) and ($iv$) are also best possible.

\noindent{\bf Definition 5.} Let $k, t, r$ be integers such that $t>r>2$, $k=(^{t-1}_{r-1})$ and $kr\geq 2t$. Assume $n$ and $l$ are integers satisfying $n=l+pt$ and $l\geq 2t+2$, where $p\geq0$. Let $a=\lceil\frac{l}{2}\rceil$ and $b=\lfloor\frac{l}{2}\rfloor$. Then $a,b\geq t+1$.

Let $H_0$ be a $r$-uniform hypergraph obtained from the disjoint union of $K_a^r$  and $K_b^r$ by adding $k$ edges joining $K_a^r$  and $K_b^r$.
Let $H$ be a star-like $r$-uniform hypergraph with the nucleus $H_0$ and $p$ $K_t^r$-satellites, adding $k$ edges from each satellite to the nucleus such that ($i$) each vertex in the satellite adjacent to some added edge (we can do this by $kr\geq 2t$); ($ii$) not all of the added $k$ edges are incident with the same complete subhypergraph  $K_a^r$  or $K_b^r$.

\begin{thm}
Let $H$ be a $r$-uniform hypergraph constructed in Definition 5. Then $H$ is $(k,l)$-edge-maximal.
\end{thm}

\noindent{\bf Proof.} By definition, there is no subhypergraph $H'$ in $H$ such that $|V(H')|\geq l$ and $\kappa'(H')\geq k+1$. We will prove the theorem by induction on $p$. If $p=0$, then $|V(H)|=l$. Since $\kappa'(K_a^r)=(^{a-1}_{r-1})\geq(^{t}_{r-1})>k$ and $\kappa'(K_b^r)=(^{b-1}_{r-1})\geq(^{t}_{r-1})>k$, the only edge-cut with $k$ edges of $H$ is these edges connecting $K_a^r$  and $K_b^r$. For any $e\in E(H^c)$, we have $e\in E_{H^c}[V(K_a^r),V(K_b^r)]$. Thus every edge-cut of $H+e$ has cardinality at least $k+1$, that is, $\kappa'(H+e)\geq k+1$. Thus $H$ is $(k,l)$-edge-maximal.

Now suppose $p\geq 1$. We assume that each hypergraph constructed in Definition 5 with less than $l+pt$ vertices is $(k,l)$-edge-maximal. In the following, we will show that each $H$ in Definition 5 with $l+pt$ vertices is also $(k,l)$-edge-maximal.

On the contrary, assume that there is an edge $e\in E(H^c)$ such that $H+e$ contains no subhypergraph $H'$ such that $|V(H')|\geq l$ and $\kappa'(H')\geq k+1$. Let $F$ be an edge-cut in $H+e$ with cardinality at most $k$.
By Definition 5 ($i$), we have $\delta(H)\geq k+1$. By (1) in page 4 and Definition 5 ($ii$), we obtain that  edge-cuts in $H+e$ with cardinality at most $k$ are these $k$ edges joining one satellite to the nucleus.
Thus there is a component, say $H_1$, of $H-F$, such that $H_1$ is the hypergraph obtained from $H$ by deleting one satellite and $e\in H_1^c$.
By induction assumption, $H_1+e$ contains a subhypergraph $H_1'$ such that $|V(H_1')|\geq l$ and $\kappa'(H_1')>k$. But $H_1'$ is also a subhypergraph of $H+e$, a contradiction.
$\Box$

\end{document}